\documentclass{amsart}

\usepackage[colorlinks=true, urlcolor=black, citecolor=black, linkcolor=black, hyperfootnotes=true]{hyperref}
\usepackage{amssymb}
\usepackage{aliascnt}
\usepackage{mathscinet}
\usepackage{enumitem}

\newcounter{dummy}
\makeatletter
\newcommand\myitem[1][]{\item[#1]\refstepcounter{dummy}\def\@currentlabel{#1}}
\makeatother

\numberwithin{equation}{section}

\newtheorem{thm}{Theorem}[section]

\newaliascnt{prp}{thm}
\newtheorem{prp}[prp]{Proposition}
\aliascntresetthe{prp}

\newaliascnt{cor}{thm}
\newtheorem{cor}[cor]{Corollary}
\aliascntresetthe{cor}

\theoremstyle{definition}

\newaliascnt{dfn}{thm}
\newtheorem{dfn}[dfn]{Definition}
\aliascntresetthe{dfn}

\newaliascnt{xpl}{thm}

\aliascntresetthe{xpl}

\newaliascnt{qst}{thm}
\newtheorem{qst}[qst]{Question}
\aliascntresetthe{qst}

\newaliascnt{rmk}{thm}
\newtheorem{rmk}[rmk]{Remark}
\aliascntresetthe{rmk}

\author{Tristan Bice}
\author{Wies\l aw Kubi\'s}
\email{tristan.bice@gmail.com}
\email{kubis@math.cas.cz}
\thanks{The authors are supported by the GA\v{C}R project EXPRO 20-31529X and RVO: 67985840 at the Institute of Mathematics of the Czech Academy of Sciences in Prague, Czech Republic}

\title{Wallman Duality for Semilattice Subbases}

\begin{document}

\begin{abstract}
We extend Wallman's classic duality from lattice bases to semilattice subbases and from compact to locally closed compact spaces.  Moreover, we make this duality functorial via appropriate relational morphisms.
\end{abstract}

\maketitle

\section*{Introduction}

\subsection*{Motivation}

This paper is an extension of the ideas in \cite{Wallman1938}.  Despite being over 80 years old, there is still much inspiration to be drawn from \cite{Wallman1938}, which historically has been somewhat overshadowed by Stone's work from around the same time.  However, their motivations were really opposite in that Stone wanted topological representations of order structures, while Wallman was after order theoretic representations of topological spaces.

From the topological point of view, Wallman's results are more appealing as they apply to quite general (e.g. connected) compact spaces commonly found in analysis.  The modern approach to point-free topology applies to general spaces too (see \cite{PicadoPultr2012}), but at the cost of working with big lattices, namely frames representing the entire open set lattice.  In \cite{Wallman1938}, Wallman showed that it actually suffices to deal with a lattice representing a basis, or even an abstract simplicial complex representing a mere subbasis, at least when dealing with compact $T_1$ spaces.

The first question that naturally arises is whether the lattice aspect can also be generalised from bases to subbases.  Somewhat surprisingly, we find this is indeed possible by in some sense interpolating between the two parts of Wallman's paper.

The next question is whether Wallman duality admits a local extension analogous to the well known locally compact generalisation of Stone duality.  Wallman's approach via closed sets somewhat obscures the potential for doing this.  However, upon translation to open sets, local generalisations become more apparent, even for weak notions of local compactness.

The next natural task is to make this duality functorial.  The counterparts of continuous functions are not semilattice homomorphisms, as one might expect from the functorial aspect of Stone duality, but rather relations between $\vee$-semilattices satisfying certain key properties related to continuity.

\subsection*{Outline}

To start with in \autoref{Grills}, we examine proper minimal non-empty grills.  These play the same vital role in $\vee$-semilattices that prime filters play in distributive lattices.  In \autoref{TheSpectrum} we then examine the spectrum of such grills and show how they recover any locally relatively compact $T_1$ space from a $\cup$-subbasis in \autoref{Recovery}.  We then make a brief detour in \autoref{Products} to show that products of bounded $\vee$-semilattices correspond nicely to products of topological spaces.

In \autoref{NearFar}, we return to the general theory by examining near and far subsets.  These play the same role in $\vee$-semilattices that subsets with or without infimum $0$ play in distributive lattices.  In particular, they allow us to extend the usual `rather below' relation $\prec$ to $\vee$-semilattices in \autoref{RatherBelow}.  We next provide a connection to Wallman's original work in \autoref{Bases} by characterising $\cup$-bases among $\cup$-subbases in an order theoretic way.  This also yields a first order characterisation of $\prec$ in \eqref{Rather1}, and a version of distributivity for $\prec$ in \autoref{precDistributive}.

The final piece of the puzzle is an extension of subfitness to (even unbounded) $\vee$-semilattices, which we examine in \autoref{Subfitness}.  This allows us to characterise $\leq$ and $\prec$ in terms of containment and closed containment in the spectrum, as seen in \autoref{Faithful} and \autoref{precCharProp}.  Moreover, subfitness allows us to show that the spectrum is locally closed compact, thus yielding a duality of relatively compact $T_1$ $\cup$-subbases with subfit round $\vee$-semilattices.  In \autoref{Functoriality} we show how to make this duality functorial by obtaining (very) continuous functions from `$\vee$-relations'.  Lastly, in \autoref{CoverRelations}, we examine cover relations as an alternative to $\vee$-semilattices if one wants to deal with arbitrary subbases rather than $\cup$-subbases.

\section{Grills}\label{Grills}

In point-free topology, prime filters play a fundamental role.  To work with subbases rather than bases, we need to consider more general grills.

\begin{dfn}
In any poset, we define \emph{grills}, \emph{ideals} and \emph{filters} as follows.
\begin{align*}
\label{Grill}\tag{Grill}p\in G\text{ or }q\in G\qquad&\Leftrightarrow\qquad\forall r\geq p,q\ (r\in G).\\
\label{Ideal}\tag{Ideal}p,q\in I\qquad&\Leftrightarrow\qquad\exists r\in I\ (r\geq p,q).\\
\label{Filter}\tag{Filter}p,q\in F\qquad&\Leftrightarrow\qquad\exists r\in F\ (r\leq p,q).
\end{align*}
\end{dfn}

Note a complement of an ideal is a grill and vice versa, and the down-directed grills are precisely the prime filters (`prime up-set' might thus be an appropriate term, although we stick to the grill terminology introduced in \cite{Choquet1947}).

We are particularly interested in $\vee$-semilattices, i.e. non-empty partially ordered sets in which every $p$ and $q$ has a supremum $p\vee q$.

\begin{center}
\textbf{From now on we assume $S$ is a $\vee$-semilattice unless otherwise stated.}
\end{center}

In this case, $G\subseteq S$ is a grill if and only if
\[p\vee q\geq g\in G\qquad\Rightarrow\qquad p\in G\text{ or }q\in G.\]
One immediately sees that this then extends to non-empty finite $F\subseteq S$, i.e.
\[\bigvee F\geq g\in G\qquad\Rightarrow\qquad F\cap G\neq\emptyset.\]

We are particularly interested in proper minimal non-empty grills.  The following first order characterisation of minimality will be used repeatedly.

\begin{prp}\label{MinimalGrillEquivalent}
If $G$ is a non-empty grill in $S$ then
\[G\text{ is proper and minimal}\qquad\Leftrightarrow\qquad\forall s\in S\ \forall g\in G\ \exists t\in S\setminus G\ (s\leq t\vee g).\]
\end{prp}

\begin{proof}
If $G$ were not minimal then we would have another non-empty grill $H\subsetneqq G$.  Taking any $s\in H$ and $g\in G\setminus H\subseteq S\setminus H$ we see that there could not be any $t\in S\setminus G\subseteq S\setminus H$ with $s\leq t\vee g$ because $H$ is a grill.  If $G$ were not proper then there could not be any $t\in S\setminus G$ whatsoever, thus proving $\Leftarrow$.

Conversely, say we had $s\in S$, $g\in G$ and $s\nleq t\vee g$, for all $t\in S\setminus G$.  As long as $G$ is proper, $S\setminus G$ and $g$ would generate an ideal $I=\bigcup_{t\in S\setminus G}(t\vee g)^\geq$ containing $g$ and avoiding $s$ and hence $S\setminus I$ would be a grill containing $s$ and avoiding $g$.  Thus $\emptyset\neq S\setminus I\subsetneqq G$ so $G$ would not minimal, thus proving $\Rightarrow$.
\end{proof}

Recall that $S$ is said to be \emph{distributive} if, for all $p,q\in S$,
\[\tag{Distributive}s\leq p\vee q\qquad\Rightarrow\qquad\exists p'\leq p\ \exists q'\leq q\ (s=p'\vee q').\]
Distributivity provides a connection between grills and filters.

\begin{prp}\label{Grill=>Filter}
If $S$ is distributive, every minimal non-empty grill $G$ is a filter.
\end{prp}

\begin{proof}
If $S$ consists of a single element then $S$ itself is the only non-empty grill, which is certainly a filter.  Otherwise, every minimal non-empty grill is necessarily proper \textendash\, any $p\in S$ which is not a maximum generates a proper ideal $p^\geq$ and hence $S\setminus p^\geq$ is a proper non-empty grill, in particular $S$ itself is not a minimal non-empty grill.  For any $g,h\in G$, \autoref{MinimalGrillEquivalent} then yields $s\in S\setminus G$ with $h\leq s\vee g$.  By distributivity, we have $s'\leq s$ and $g'\leq g$ with $h=s'\vee g'$.  In particular, $g'\leq g,h$.  Note $s'\notin G$, as $s'\leq s\notin G$, and hence $g'\in G$, as $s'\vee g'=h\in G$.
\end{proof}

\section{The Spectrum}\label{TheSpectrum}

\begin{dfn}
The \emph{spectrum} of $S$ is given by
\[\widehat{S}=\{G\subseteq S:G\text{ is a proper minimal non-empty grill}\}.\]
\end{dfn}

As noted in the proof of \autoref{Grill=>Filter}, `proper' here is superfluous as long as $S$ contains at least two distinct elements.

We consider $\widehat{S}$ as a space with the topology generated by $(\widehat{S}_g)_{g\in S}$ where
\[\widehat{S}_g=\{G\in\widehat{S}:g\in G\},\]
i.e. we are taking $(\widehat{S}_g)_{g\in S}$ as a subbasis for the topology on $\widehat{S}$.  Our goal in this section is to show that the spectrum allows us to recover a large class of spaces from the semilattice structure of a subbasis, at least among $T_1$ spaces.

Let us call a family $\mathcal{P}\subseteq\mathcal{P}(X)$ of subsets \emph{$T_1$} if, for all singleton or empty $a,b\subseteq X$,
\[\tag{$T_1$}a\cap b=\emptyset\qquad\Rightarrow\qquad\exists p\in\mathcal{P}\ (a\subseteq p\text{ and }b\cap p=\emptyset).\]
If $X$ contains at least two points then it suffices to consider singleton $a$ and $b$, but if $X$ itself is a singleton or empty then taking $a$ and/or $b$ to be empty shows that
\begin{align*}
X=\{x\}\qquad&\Rightarrow\qquad\{\emptyset,X\}\text{ is the only $T_1$ family}.\\
X=\emptyset\qquad&\Rightarrow\qquad\{\emptyset\}=\{X\}\text{ is the only $T_1$ family}.
\end{align*}
A space $X$ is $T_1$ iff it has a $T_1$ subbasis, and if $\mathcal{P}$ is a subbasis of a $T_1$ space $X$ then
\[\mathcal{P}\text{ is }T_1\qquad\Leftrightarrow\qquad\bigcup\mathcal{P}=X\text{ and }\bigcap\mathcal{P}=\emptyset.\]
Again, this is automatic if $X$ contains at least two points.

\begin{prp}\label{MinimalT1}
The spectrum $\widehat{S}$ is a $T_1$ space with $T_1$ subbasis $(\widehat{S}_p)_{p\in S}$.
\end{prp}

\begin{proof}
For any distinct $G,H\in\widehat{S}$, minimality yields $g\in G\setminus H$ and $h\in H\setminus G$ so $G\in\widehat{S}_g\setminus\widehat{S}_h$ and $H\in\widehat{S}_h\setminus\widehat{S}_g$, showing that $\widehat{S}$ is a $T_1$ space.  Moreover, $\widehat{S}=\bigcup_{p\in S}\widehat{S}_p$, as each $G\in\widehat{S}$ is non-empty, and $\emptyset=\bigcap_{p\in S}\widehat{S}_p$, as each $G\in\widehat{S}$ is proper.
\end{proof}

In fact, any $T_1$ space arises in this way as long as it satisfies a weak version of local compactness.  First, given a topology $\mathcal{O}(X)$ on a set $X$, let us call $R\subseteq X$ \emph{relatively compact} if every open cover of $X$ has a finite subcover of $R$, i.e.
\[\mathcal{C}\subseteq\mathcal{O}(X)\quad\text{and}\quad X\subseteq\bigcup\mathcal{C}\qquad\Rightarrow\qquad\exists\text{ finite }\mathcal{F}\subseteq\mathcal{C}\ (R\subseteq\bigcup\mathcal{F}).\]
Equivalently, $R$ is relatively compact if every ultrafilter in $R$ converges in $X$, which is a standard notion in convergence theory \textendash\, see \cite{Lowen-Colebunders1983} and \cite[Ch IX]{DoleckiMynard2016}.

\begin{rmk}
This is not equivalent to saying $\mathrm{cl}(R)$ is compact, which is another commonly cited notion of relative compactness.  One the one hand, if $R$ is contained in a compact subset then certainly $R$ is relatively compact, and in regular spaces or locally closed compact spaces (see below) every relatively compact set has compact closure.  However, even in Hausdorff spaces, a set can be relatively compact but not contained in any compact subset, e.g. in the half-disc topology on the upper half plane \textendash\, see \cite[Counterexample 78]{SteenSeebach1978}.
\end{rmk}

We distinguish the following weak notions of local compactness.

\begin{dfn}
We call a topological space $X$
\begin{enumerate}[leftmargin=15pt]
\item \emph{locally relatively compact} if each $x\in X$ has a relatively compact neighbourhood.
\item \emph{locally closed compact} if each $x\in X$ has a closed compact neighbourhood.
\end{enumerate}
\end{dfn}

\begin{dfn}
We call $S\subseteq\mathcal{O}(X)$ \emph{$\cup$-closed} if $S$ is closed under pairwise unions, i.e. if $S$ is a $\vee$-subsemilattice of $\mathcal{O}(X)$.  We call a $\cup$-closed subbasis a \emph{$\cup$-subbasis}.
\end{dfn}

The import of the following result is that any locally relatively compact $T_1$ space can be recovered from an appropriate subbasis, ordered by inclusion $\subseteq$.

\begin{thm}\label{Recovery}
For any relatively compact $T_1$ $\cup$-subbasis $S$ of a space $X$,
\[x\mapsto S_x=\{s\in S:x\in s\}\]
is a homemorphism from $X$ onto $\widehat{S}$.
\end{thm}

\begin{proof}
As $S$ is $T_1$, each $S_x$ is a proper non-empty grill.  If it were not minimal then we would have another non-empty grill $G\subsetneqq S_x$.  Again as $S$ is $T_1$, for each $y\neq x$, we have $s_y\in S$ with $y\in s_y\not\ni x$.  Taking any $g\in G\subseteq S_x$ and $s\in S_x\setminus G$, we see that $X=s\cup\bigcup_{y\neq x}s_y$ and hence $g\subseteq s\cup\bigcup_{y\in F}s_y$, for some finite $F\subseteq X\setminus\{x\}$, as $g$ is relatively compact.  But as $s\notin G$ and each $s_y\notin S_x\supseteq G$, this contradicts the fact $G$ is a grill.  Thus $S_x$ is minimal.

Conversely, to show there are no other minimal non-empty grills, it suffices to show every non-empty grill $G$ contains $S_x$, for some $x\in X$.  If this were not the case then $S\setminus G$ would cover $X$.  Taking any $g\in G$, relative compactness would then yield finite $F\subseteq S\setminus G$ with $g\subseteq\bigcup F$, again contradicting the fact $G$ is a grill.

Thus $x\mapsto S_x$ is a bijection from $X$ to $\widehat{S}$.  Consequently, it is a homeomorphism, as we immediately see that it maps the subbasis $S$ onto the subbasis $(\widehat{S}_p)_{p\in S}$.
\end{proof}

When $X$ above is compact and $S$ is also closed under pairwise intersections (and is hence a basis), the resulting posets can be characterised as the subfit bounded distributive lattices, as shown in \cite{Wallman1938}.  The ultimate goal would be to extend this to the more general $\vee$-semilattices above.

\begin{qst}\label{RelativeQuestion}
Is there an order theoretic characterisation of $\vee$-semilattices arising from relatively compact $T_1$ $\cup$-subbases?
\end{qst}

If we restrict our attention slightly to locally closed compact spaces then the answer is yes, specifically the round subfit $\vee$-semilattices are precisely those arising from relatively compact $T_1$ $\cup$-subbases in locally closed compact spaces.  Proving this (and extending `subfit' and `round' to $\vee$-semilattices) is our primary goal.

\section{Products}\label{Products}
Before moving on, however, let us point out that the spectrum behaves well with respect to products, at least for bounded $\vee$-semilattices.  Indeed, this is one advantage of semilattice subbases over lattice bases (there is a corresponding coproduct of frames, but it is somewhat more involved \textendash\, see \cite[Ch IV]{PicadoPultr2012}).

For motivation, say we have $\cup$-subbases $S$ and $S'$ of spaces $X$ and $X'$ and note
\[\{(s\times X')\cup(X\times s'):s\in S\text{ and }s'\in S'\}\]
then forms a $\cup$-subbasis of $X\times X'$, as long as $\emptyset\in S$ and $\emptyset\in S'$.  Indeed,
\[((s\times X')\cup(X\times s'))\cup((t\times X')\cup(X\times t'))=((s\cup t)\times X')\cup(X\times(s'\cup t')).\]
Also note that the sets $((s\times X')\cup(X\times s'))_{s\in S,s'\in X'}$ are all distinct except when $s=X$ or $s'=X'$, in which case we get $X\times X'$, i.e.
\[(s\times X')\cup(X\times X')=X\times X'=(X\times X')\cup(X\times s').\]

Accordingly, given bounded $\vee$-semilattices $S$ and $S'$ with minima $0$ and $0'$ and maxima $1$ and $1'$ respectively, let $[S\times S']$ denote the usual product $S\times S'$ where we identify pairs containing $1$ or $1'$.  More precisely, let
\[1^\times=(\{1\}\times S')\cup(S\times\{1'\})=\{(s,s'):s=1\text{ or }s'=1'\}\]
and let $[S\times S']=\{[s,s']:s\in S\text{ and }s'\in S'\}$ where
\[[s,s']=\begin{cases}\{(s,s')\}&\text{if }(s,s')\notin1^\times.\\ 1^\times&\text{if }(s,s')\in1^\times.\end{cases}\]
Then $[S\times S']$ is again a $\vee$-semilattice where
\[[s,s']\vee[t,t']=[s\vee t,s'\vee t'].\]

\begin{thm}
For any bounded $\vee$-semilattices $S$ and $S'$, the map
\[(G,G')\mapsto G\sqcup G'=[G\times S']\cup[S\times G']=\{[s,s']:s\in G\text{ or }s'\in G'\}\]
is a homeomorphism from $\widehat{S}\times\widehat{S'}$ onto $\widehat{[S\times S']}$.
\end{thm}

\begin{proof}
Take $G\in\widehat{S}$ and $G'\in\widehat{S'}$.  Certainly $G\sqcup G'$ is an up-set.  If $[p\vee q,p'\vee q']=[p,p']\vee[q,q']\in G\sqcup G'$ then either $p\vee q\in G$ or $p'\vee q'\in G'$.  In the former case, either $p\in G$ or $q\in G$, and hence either $[p,p']\in[G\times S']$ or $[q,q']\in[G\times S]$.  Likewise, in the latter case either $[p,p']\in[S\times G']$ or $[q,q']\in[S\times G']$.  This shows that $G\sqcup G'$ is a grill and it is also non-empty and proper, as both $G$ and $G'$ are.  Now if $[g,p']\in G\times S'$ then \autoref{MinimalGrillEquivalent} yields $t\in S\setminus G$ with $1=t\vee g$ and hence $1^\times=[t\vee g,p']=[t,0']\vee[g,p']$.  Note $[t,0']\notin G\sqcup G'$, as $t\notin G$ and $0'\notin G'$.  Likewise, if $[p,g']\in S\times G'$ then we have $[0,t']\in[S\times S']\setminus G\sqcup G'$ with $1^\times=[p,g']\vee[0,t']$ and hence $G\sqcup G'$ is minimal, by \autoref{MinimalGrillEquivalent}, i.e. $G\sqcup G'\in\widehat{[S\times S']}$.

Conversely, take any $G^\times\in\widehat{[S\times S']}$ and let
\[G=\{g\in S:[g,0']\in G^\times\}\qquad\text{and}\qquad G'=\{g'\in S':[0,g']\in G^\times\}.\]
We immediately see that $G$ and $G'$ are proper non-empty grills.  For minimality, say $g\in G$ so $[g,0']\in G^\times$ and hence we have $[t,t']\in[S\times S']\setminus G^\times$ such that $1^\times=[g,0']\vee[t,t']=[g\vee t,t']$.  As $[t,t']\notin G^\times$, $1^\times\neq[t,t']$ and hence $t'\neq1'$.  So the only way we could have $[g\vee t,t']=1^\times$ is if $g\vee t=1$.  Also note $[t,0']\notin G^\times$, as $[t,0']\leq[t,t']\notin G^\times$, and hence $t\notin G$.  By \autoref{MinimalGrillEquivalent}, this shows that $G$ is minimal and, likewise, $G'\in\widehat{S'}$.  As $[G\times\{0'\}]\cup[\{0\}\times G']$ is coinitial in $G\sqcup G'=[G\times S']=[S\times G']$ and $G^\times$ is an up-set, it follows that $G\sqcup G'\subseteq G^\times$.  On the other hand, if $s\in S\setminus G$ and $s'\in S'\setminus G'$, i.e. $[s,0'],[0,s']\notin G^\times$, then $[s,s']=[s,0']\vee[0,s']\notin G^\times$, as $G^\times$ is a grill.  Thus $G^\times=G\sqcup G'$, showing that $(G,G')\mapsto G\sqcup G'$ takes $\widehat{S}\times\widehat{S'}$ onto $\widehat{[S\times S']}$.  This argument also shows that $(G,G')\mapsto G\sqcup G'$ maps $O_p\times\widehat{S'}$ and $\widehat{S}\times O_{p'}$ onto $O_{[p,0']}$ and $O_{[0,p']}$ respectively.  As these are subbases of $\widehat{S}\times\widehat{S'}$ and $\widehat{[S\times S']}$, the map is also a homeomorphism.
\end{proof}

\section{Near vs Far}\label{NearFar}

\begin{dfn}
We call any non-empty finite $F\subseteq S$ \emph{near} or \emph{far} if
\begin{align*}
\label{Near}\tag{Near}&\exists p,q\in S\ (\forall f\in F\ (p\leq q\vee f)\text{ and }p\nleq q).\\
\label{Far}\tag{Far}&\forall p,q\in S\ (\forall f\in F\ (p\leq q\vee f)\ \Rightarrow\ p\leq q).
\end{align*}
\end{dfn}

As the terminology suggests, these are opposite concepts, i.e. near means not far and vice versa.  Intuitively, $F$ is near if it has `non-empty intersection'.

\begin{prp}\label{SubbasisNear}
If $S\supseteq F$ is a relatively compact $T_1$ $\cup$-subbasis of $X$ then
\[F\text{ is near}\qquad\Leftrightarrow\qquad\bigcap F\neq\emptyset.\]
\end{prp}

\begin{proof}
If $\bigcap F=\emptyset$ then $F$ is far, as $p\subseteq q\cup f$, for all $f\in F$, implies
\[p\subseteq\bigcap_{f\in F}(q\cup f)=q\cup\bigcap F=q\cup\emptyset=q.\]

Conversely, say we have $x\in\bigcap F$.  As $S$ is $T_1$, we have $p\in S$ containing $x$ and, for every $y\notin\bigcap F\ni x$, we have $q\in S$ with $y\in q\not\ni x$.  As $p$ is relatively compact, finitely many such $q$ cover $p\setminus\bigcap F$.  As $S$ is $\cup$-closed, we can take their union to obtain a single such $q$.  Thus $p\subseteq(\bigcap F)\cup q\subseteq f\cup q$, for all $f\in F$, even though $p\nsubseteq q$, as $x\in p\setminus q$, i.e. $F$ is near.
\end{proof}

We have a very similar characterisation for distributive $\vee$-semilattices.

\begin{prp}\label{DistributiveNear}
If $S$ is distributive with minimum $0$ then, for any finite $F\subseteq S$,
\[F\text{ is near}\qquad\Leftrightarrow\qquad\bigwedge F\neq0.\]
\end{prp}

\begin{proof}
If $\bigwedge F\neq0$ then we have some non-zero $f'\leq f$, for all $f\in F$, so we can take $p=f'$ and $q=0$ in the definition above to show that $F$ is near.

Conversely, if $F$ is near then we have $p,q\in S$ with $q\ngeq p\leq q\vee f$, for all $f\in F$.  Taking any $f_1\in F$, distributivity yields $q'\leq q$ and $p_1\leq f_1$ with $p=q'\vee p_1$.  As $p\nleq q$, it follows that $p_1\nleq q$, however $p_1\leq p\leq q\vee f$, for all $f\in F\setminus\{f_1\}$.  Continuing in this way, we obtain $p_n\leq f$, for all $f\in F$, with $p_n\nleq q$ so $\bigwedge F\neq0$.
\end{proof}

Finite subsets of proper minimal non-empty grills are always near.

\begin{prp}\label{NearMinimalGrill}
If $G\in\widehat{S}$ then every finite $F\subseteq G$ is near.
\end{prp}

\begin{proof}
Take any $g\in G$.  For all $f\in F\subseteq G$, \autoref{MinimalGrillEquivalent} yields $f'\in S\setminus G$ with $g\leq f\vee f'$.  Let $t=\bigvee_{f\in F}f'$ so $g\nleq t\in S\setminus G$, as $G$ is a grill, and $g\leq f\vee f'\leq f\vee t$, for all $f\in F$, showing that $F$ is indeed near.
\end{proof}

\begin{prp}\label{pGrill}
For any $p\in S$ and any grill $G$ containing $p$, we have another grill $H\subseteq G$ containing $p$ such that $F$ is near, for all finite $F\subseteq H$.
\end{prp}

\begin{proof}
By Kuratowski-Zorn, we have a grill $H\subseteq G$ that is minimal among grills containing $p$.  Say we had finite $F\subseteq H$ with $F$ far.  For all $q\in S\setminus H$, we have $p\nleq q$, as $p\in H$, so the definition of far yields $f\in F$ such that $p\nleq q\vee f$.  In fact, we claim that some $f\in F$ must satisfy $p\nleq q\vee f$, for all $q\in S\setminus H$.  If not, for every $f\in F$, we would have $f'\in S\setminus H$ with $p\leq f'\vee f$.  Taking $F'=\{f':f\in F\}$ would then yield $p\leq\bigvee F'\vee f$, for all $f\in F$, even though $\bigvee F'\in S\setminus H$, as $H$ is a grill, contradicting our earlier observation.  Now the claim is proved, taking $f\in F$ such that $p\nleq q\vee f$, for all $q\in S\setminus H$, we see that the ideal $I$ generated by $f$ and $S\setminus H$ avoids $p$ and hence $S\setminus I$ is a grill containing $p$ and avoiding $f$.  Thus $S\setminus I\subsetneqq H$, as $f\in F\subseteq H$, contradicting minimality.  Thus $F$ must have been near.
\end{proof}

Just finally, we make a couple of elementary observations about far subsets.  Let
\begin{align*}
E\vee F&=\{e\vee f:e\in E\text{ and }f\in F\}.\\
E\geq F\quad&\Leftrightarrow\quad\forall e\in E\ \exists f\in F\ (e\geq f).
\end{align*}

\begin{prp}
For any finite $E,F\subseteq S$,
\begin{align}
\label{veeFar}E\text{ and }F\text{ are far}\qquad&\Rightarrow\qquad E\vee F\text{ is far}.\\
\label{geqFar}E\text{ is far and }E\geq F\qquad&\Rightarrow\qquad F\text{ is far}.
\end{align}
\end{prp}

\begin{proof}\
\begin{itemize}
\item[\eqref{veeFar}] If $E$ and $F$ are far and $p\leq q\vee e\vee f$, for all $e\in E$ and $f\in F$, then $p\leq q\vee e$, for all $e\in E$, as $F$ is far, and hence $p\leq q$, as $E$ is far.

\item[\eqref{geqFar}] If $E$ is far, $E\geq F$ and $p\leq q\vee f$, for all $e\in E$, then $p\leq q\vee e$, for all $e\in E$, and hence $p\leq q$, as $E$ is far. \qedhere
\end{itemize}
\end{proof}

\begin{prp}
$0$ is a minimum of $S$ iff $\{0\}$ is far.
\end{prp}

\begin{proof}
If $0$ is a minimum of $S$ then, for any $p,q\in S$, $p\leq q\vee 0$ implies $p\leq q$, as $q\vee 0=q$, showing that $\{0\}$ is far.  Conversely, if $f$ is not a minimum of $S$ then we have some $p\ngeq f$.  As $f\leq p\vee f$, this shows that $\{f\}$ is near.
\end{proof}

\section{Rather Below}\label{RatherBelow}

\begin{dfn}\label{RatherBelowDefinition}
We define the \emph{rather below} relation $\prec$ on $S$ by
\[p\prec q\qquad\Leftrightarrow\qquad\forall s\nleq q\ \exists\text{ finite }F\subseteq S\ (F\cup\{p\}\text{ is far and }\forall f\in F\ (s\leq f\vee q)).\]
\end{dfn}

Note $s$ can always be added to $F$ above, so it suffices to consider non-empty $F$.

\subsection{Equivalents}
In certain situations, $\prec$ has various equivalent characterisations.

\begin{prp}\label{precMax}
If $S$ has a maximum $1$, i.e. $p\leq1$, for all $p\in S$, then
\begin{align*}
p\prec q\qquad&\Leftrightarrow\qquad q=1\text{ or }\exists\text{ finite }F\subseteq S\ (F\cup\{p\}\text{ is far and }\forall f\in F\ (1\leq f\vee q)).\\
1\prec q\qquad&\Leftrightarrow\qquad q=1.
\end{align*}
\end{prp}

\begin{proof}
If $q\neq1$ then it suffices to take $s=1$ in the formula defining $p\prec q$, as $1\leq f\vee q$ implies $s\leq f\vee q$, for all $s\in S$.  On the other hand, the formula defining $p\prec1$ holds vacuously for all $s\nleq1$.

The first $\Leftarrow$ immediately yields the second $\Leftarrow$.  Conversely, say $1\prec q\neq1$.  Then we can take $s=1$ in the formula defining $1\prec q$, which means we have finite $F\subseteq S$ such that $F\cup\{1\}$ is far and $1\leq f\vee q$, for all $f\in F$.  But this means $F$ is also far and hence $1\leq q$, a contradiction, proving $\Rightarrow$.
\end{proof}

Note the last equivalence is saying `every $\leq$-maximum is a unique $\prec$-maximum'.

Usually, we can also replace $\forall s\nleq q$ with $\forall s\in S$ in the definition of $\prec$.

\begin{prp}\label{foralls}
If $S$ has no maximum or has at least one far subset then
\[p\prec q\qquad\Leftrightarrow\qquad\forall s\in S\ \exists\text{ finite }F\subseteq S\ (F\cup\{p\}\text{ is far and }\forall f\in F\ (s\leq f\vee q)).\]
\end{prp}

\begin{proof}
If $p,q,s\in S$ and $q$ is not a maximum then we have $t\nleq q$ and hence $s\vee t\nleq q$.  Thus $p\prec q$ implies the existence of finite $F\subseteq S$ with $F\cup\{p\}$ far and $s\leq s\vee t\leq f\vee q$, for all $f\in F$.  On the other hand, if $q=1$ is the maximum of $S$ and $S$ has a far subset $F$ then certainly $s\leq f\vee1$, for all $f\in F$.
\end{proof}

However, there are $S$ for which the above result does not apply, e.g. $S=\mathbb{N}$ in the reverse ordering or $S=$ the cofinite subsets of an infinite set ordered by inclusion.

In distributive $\vee$-semilattices, we can replace finite subsets with singletons.

\begin{prp}
If $S$ is distributive with minimum $0$ then
\[p\prec q\qquad\Leftrightarrow\qquad\forall s\in S\ \exists t\in S\ (p\wedge t=0\text{ and }s\leq t\vee q).\]
\end{prp}

\begin{proof}
If the right hand side holds then $p\prec q$, as witnessed by taking $F=\{t\}$ in the definition of $\prec$, noting that $p\wedge t=0$ means $\{p,t\}$ is far, by \autoref{DistributiveNear}.

Conversely, say $p\prec q$.  If $q=1$ is a maximum of $S$ then the right hand side holds with $t=0$.  Otherwise, for any $s\in S$, we have finite $F\subseteq S$ such that $F\cup\{p\}$ is far and $s\leq f\vee q$, for all $f\in F$.  Taking $f_1\in F$, distributivity yields $s_1\leq f_1$ and $q_1\leq q$ with $s=s_1\vee q_1$.  Taking any other $f_2\in F$, we see that $s_1\leq s\leq f_2\vee q$ so distributivity again yields $s_2\leq f_2$ and $q_2\leq q$ with $s_1=s_2\vee q_2$ and hence $s=s_1\vee q_1=s_2\vee q_2\vee q_1$.  Continuing in this way we obtain $s_n\leq f$, for all $f\in F$, such that $s\leq s_n\vee q$.  As $F\cup\{p\}$ is far, we must have $p\wedge s_n\leq p\wedge\bigwedge F=0$, by \autoref{DistributiveNear}, i.e. we can take $t=s_n$ above.
\end{proof}

\begin{rmk}
If $S$ also has a maximum $1$ then it suffices to take $s=1$ above, i.e.
\[p\prec q\qquad\Leftrightarrow\qquad\exists t\in S\ (p\wedge t=0\text{ and }1\leq t\vee q),\]
so $\prec$ agrees with the rather below relation defined for frames in \cite[V.5.2]{PicadoPultr2012}, originally called `well-inside' and defined for distributive lattices in \cite[III.1.1]{Johnstone1986}.
\end{rmk}

\subsection{Properties}
Now we examine some properties of $\prec$, the first being auxiliarity.

\begin{prp}
For any $p,p',q,q'\in S$,
\[\tag{Auxiliarity}\label{Auxiliarity}p\leq p'\prec q'\leq q\qquad\Rightarrow\qquad p\prec q\qquad\Rightarrow\qquad p\leq q.\]
\end{prp}

\begin{proof}
Say $p\leq p'\prec q'$.  If $p'=1$ then $q'=1$, as noted above.  Otherwise, for all $s\in S$, we have finite $F\subseteq S$ with $F\cup\{p'\}$ and $p\leq p'\leq f\vee q$, for all $f\in F$.  As $p\leq p'$, $F\cup\{p\}$ is also far so this shows that $p\prec q'$.  On the other hand, if $p'\prec q'\leq q$ then certainly $p'\prec q$, as $p'\leq f\vee q'$ implies $p'\leq f\vee q$, for any $f\in S$.

Now say $p\prec q$.  If $q=1$ then certainly $p\leq q$.  Otherwise, we have $r\nleq q$ and we can take $s=p\vee r$ to get finite $F\subseteq S$ such that $F\cup\{p\}$ is far and $p\leq s\leq f\vee q$, for all $f\in F$.  Certainly $p\leq p\vee q$ too so, as $F\cup\{p\}$ is far, $p\leq q$.
\end{proof}

\begin{prp}\label{FarBelow}
If $F\cup\{p\}$ is far then, for all $q\in S$,
\[\forall f\in F\ (p\prec f\vee q)\qquad\Rightarrow\qquad p\prec q.\]
\end{prp}

\begin{proof}
As $S$ has a far subset, we can use the equivalent of $\prec$ in \autoref{foralls}.  So, for any $f\in F$ and $s\in S$, we have finite $E_f\subseteq S$ such that $E_f\cup\{p\}$ is far and $s\leq e\vee f\vee q$, for all $e\in E_f$.  Enumerate $F=\{f_1,\cdots,f_n\}$ and let
\[D=E_{f_1}\vee\cdots\vee E_{f_n}\vee F\]
so $s\leq d\vee q$, for all $d\in D$, and $D\cup\{p\}$ is far, by \eqref{veeFar} and \eqref{geqFar}.  Thus $p\prec q$.
\end{proof}

This allows us to prove an analog of \autoref{pGrill}.

\begin{prp}\label{pqGrill}
Whenever $p\nprec q$ and $p\prec p'$, we have $G\in\widehat{S}$ such that $p'\in G\not\ni q$ and $F\cup\{p\}$ is near, for all finite $F\subseteq G$.
\end{prp}

\begin{proof}
If $p=1$ is a maximum of $S$ then $p'=1$ too and $q\neq 1$ (see \autoref{precMax}).  Then Kuratowski-Zorn yields a minimal grill $G\subseteq S\setminus q^\geq$ containing $1=p=p'$ in which every finite subset is near, by \autoref{NearMinimalGrill}.

Otherwise, say $p\nprec q$ and let $I$ be an ideal containing $q$ such that $p\nprec r$, for all $r\in I$, which is also maximal with respect to this property (again using Kuratowski-Zorn).  Let $G=S\setminus I$, noting that $p'\in G$ as $p'\notin I$, and take any finite $F\subseteq G$.  By the maximality property of $I$, for every $f\in F$, we have $r_f\in I$ such that $p\prec r_f\vee f$ and hence $p\prec r\vee f$, where $r=\bigvee_{f\in F}r_f\in I$.  If $F\cup\{p\}$ were far, \autoref{FarBelow} would yield $p\prec r$, contradicting $r\in I$.  Thus $F\cup\{p\}$ is near, for all finite $F\subseteq G$.

If $G$ were not minimal, we could take minimal $H\subsetneqq G$ and $g\in G\setminus H$.  The maximality property of $I$ would then yield $r\in I$ with $p\prec g\vee r$.  Taking any $h\in H$, this means we have $F\subseteq S$ such that $F\cup\{p\}$ is far and $h\leq f\vee g\vee r$, for all $f\in F$.  Thus $F\subseteq H\subseteq G$, as $H$ is a grill containing $h$ but avoiding $g$ and $r$, contradicting what we just proved.  Thus $G$ is minimal, i.e. $G\in\widehat{S}$.
\end{proof}

\begin{prp}
For all $p,q,r\in S$,
\begin{equation}\label{pqr}
p\prec q\vee r\quad\text{and}\quad q\prec p\vee r\qquad\Leftrightarrow\qquad p\vee q\prec q\vee r,p\vee r.
\end{equation}
\end{prp}

\begin{proof}
We immediately get $\Leftarrow$ from \eqref{Auxiliarity}.

Conversely, say $p\prec q\vee r$ and $q\prec p\vee r$.  If $q\vee r=1=p\vee r$ then certainly $p\vee q\prec q\vee r,p\vee r$.  Otherwise $S$ has a far subset so we can consider the equivalent of $\prec$ in \autoref{foralls}.  Accordingly, take any $s\in S$.  As $p\prec q\vee r$, we have finite $E\subseteq S$ such that $E\cup\{p\}$ is far and $s\leq e\vee q\vee r$, for all $e\in E$.  Likewise, we have finite $F\subseteq S$ such that $F\cup\{q\}$ is far and $s\leq f\vee p\vee r$, for all $f\in F$.  Letting
\[D=E\cup\{f\vee p:f\in F\},\]
we claim that $D\cup\{p\vee q\}$ is far.  Indeed, if $a\leq b\vee p\vee q$ and $a\leq b\vee f\vee p$, for all $f\in F$, then $a\leq b\vee p$, as $F\cup\{q\}$ is far.  If $a\leq b\vee e$, for all $e\in E$ too, then $a\leq b$, as $E\cup\{p\}$ is far, proving the claim.  Now note that $s\leq e\vee q\vee r$, for all $e\in E$, and $s\leq f\vee p\vee r\leq f\vee p\vee q\vee r$, for all $f\in F$, i.e. $s\leq d\vee q\vee r$, for all $d\in D$.  This shows that $p\vee q\prec q\vee r$, while a symmetric argument yields $p\vee q\prec p\vee r$.
\end{proof}

\subsection{Roundness}

Next we show that locally closed compact spaces are closely related to `round' $\vee$-semilattices.  Answering \autoref{RelativeQuestion} for locally relatively compact spaces would thus require finding a suitable replacement for roundness.

\begin{dfn}
We call $S$ \emph{round} if, for all $p\in S$, we have $q\succ p$, for some $q\in S$.
\end{dfn}

\begin{prp}\label{T1round}
If $S$ is a $T_1$ $\cup$-subbasis of $X$ and $\mathrm{cl}(p)$ is compact, for all $p\in S$, then $S$ is round and, for all $p,q\in S$,
\[p\prec q\qquad\Leftrightarrow\qquad\mathrm{cl}(p)\subseteq q.\]
\end{prp}

\begin{proof}
If $q=X$ then certainly $p\prec q$ and $\mathrm{cl}(p)\subseteq q$, so we may assume that $q\neq X$.

Say $p\prec q$ and take $r\nleq q$.  As $\mathrm{cl}(p)$ is compact, we can cover it with finitely many subbasic sets and take $s\in S$ to be their union together with $r$.  As $p\prec q$, we have finite $F\subseteq S$ with $p\cap\bigcap F=\emptyset$ (see \autoref{SubbasisNear}) and $s\subseteq f\cup q$, for all $f\in F$.  Thus $\mathrm{cl}(p)\cap\bigcap F=\emptyset$ and $\mathrm{cl}(p)\subseteq s\subseteq\bigcap_{f\in F}(f\cup q)=(\bigcap F)\cup q$ so $\mathrm{cl}(p)\subseteq q$.

Conversely, say $\mathrm{cl}(p)\subseteq q$ so, for any $x\in X\setminus q$, we have finite $F\subseteq S$ with $x\in\bigcap F$ and $p\cap\bigcap F=\emptyset$.  For any $s\nleq q$, $\mathrm{cl}(s)\setminus q\neq\emptyset$ is compact so we have $F_1,\ldots,F_n$ with $\mathrm{cl}(s)\subseteq q\cup\bigcup_{k=1}^n\bigcap F_n$ and $p\cap\bigcap F_k=\emptyset$, for all $k\leq n$.  Letting
\[F=\big\{\bigcup_{k=1}^nf_k:\forall k\leq n\ (f_k\in F_k)\big\},\]
we have $\bigcap F=\bigcup_{k=1}^n\bigcap F_n$ so $\mathrm{cl}(s)\subseteq q\cup\bigcap F$ and $p\cap\bigcap F=\emptyset$.  Thus $F\cup\{p\}$ is far, by \autoref{SubbasisNear}, and $s\subseteq q\cup f$, for all $f\in F$, i.e. $p\prec q$.
\end{proof}

In Hausdorff spaces, $\cup$-subbases satisfy the following extra condition.

\begin{prp}
If $S$ is a $\cup$-subbasis of Hausdorff $X$ and $\mathrm{cl}(p)\subseteq q\cup r$ is compact,
\[\exists\text{ finite }F,G\subseteq S\ (\forall f\in F\ (p\subseteq f\cup r),\ \forall g\in G\ (p\subseteq q\cup g)\text{ and }\bigcap F\cap\bigcap G=\emptyset).\]
\end{prp}

\begin{proof}
As $(\mathrm{cl}(p)\setminus q)\cap(\mathrm{cl}(p)\setminus r)=\emptyset$, this follows from the well known fact that disjoint compact sets can be separated by disjoint open sets in Hausdorff spaces.
\end{proof}

Conversely, an abstract version of this implies that the spectrum is Hausdorff.

\begin{prp}
If $S$ is round and, for all $p,q,r\in S$ with $p\prec q\vee r$,
\begin{equation}\label{Hausdorff}
\exists\text{ finite }F,G\subseteq S\ (\forall f\in F\ (p\leq f\vee r),\ \forall g\in G\ (p\leq q\vee g)\text{ and }F\cup G\text{ is far}),
\end{equation}
then $\widehat{S}$ is Hausdorff.
\end{prp}

\begin{proof}
Say this condition holds and take distinct $G,H\in\widehat{S}$.  Taking the join of any elements in $G$ and $H$, we obtain $p\in G\cap H$.  Taking any $g\in G\setminus H$, \autoref{MinimalGrillEquivalent} and the fact $S$ is $\succ$-round yields $h\in H\setminus G$ with $p\prec g\vee h$.  Now take finite $E,F\subseteq S$ such that $E\cup F$ is far, $p\leq e\vee h$, for all $e\in E$, and $p\leq g\vee f$, for all $f\in F$.  Thus $E\subseteq G$, as $p\in G\not\ni h$, and $F\subseteq H$, as $p\in H\not\ni g$, i.e. $G\in\widehat{S}_E=\bigcap_{e\in E}\widehat{S}_e$ and $H\in\widehat{S}_F=\bigcap_{f\in F}\widehat{S}_f$.  But $E\cup F$ is far so $\widehat{S}_E\cap\widehat{S}_F=\emptyset$, by \autoref{NearMinimalGrill}, i.e. these open neighbourhoods of $G$ and $H$ are disjoint, showing that $\widehat{S}$ is Hausdorff.
\end{proof}

\section{Bases}\label{Bases}

Wallman's original duality concerned bases rather than subbases.  Here we show how to characterise $\cup$-bases among $\cup$-subbases using their order structure.

\begin{prp}
If $S$ is a $\cup$-basis of $X$ and $\mathrm{cl}(p\cap q)$ is compact,
\[\mathrm{cl}(p)\subseteq q\cup s\text{ and }\mathrm{cl}(q)\subseteq p\cup s\quad\Rightarrow\quad\exists r\in S\ (r\subseteq p\cap q\text{ and }p\cup q\subseteq r\cup s).\]
\end{prp}

\begin{proof}
First note
\[\mathrm{cl}(p\cap q)\subseteq\mathrm{cl}(p)\cap\mathrm{cl}(q)\subseteq(p\cap q)\cup s.\]
As $\mathrm{cl}(p\cap q)$ is compact, $\mathrm{cl}(p\cap q)\setminus s$ is a compact subset of $p\cap q$.  As $S$ is a basis, we can cover $\mathrm{cl}(p\cap q)\setminus s$ by $r\in S$ with $r\subseteq p\cap q$.  By compactness, finitely many such sets suffice.  Taking their union, we see that a single $r\in S$ with $r\subseteq p\cap q$ suffices.  As $\mathrm{cl}(p\cap q)\setminus s\subseteq r$, $p\cap q\subseteq\mathrm{cl}(p\cap q)\subseteq r\cup s$ so, as $p\subseteq\mathrm{cl}(p)\subseteq q\cup s$,
\[p\subseteq(p\cap q)\cup s\subseteq r\cup s.\]
Likewise $q\subseteq r\cup s$ so $p\cup q\subseteq r\cup s$, as required.
\end{proof}

\begin{prp}
If $S$ is round then the following conditions are equivalent.
\begin{align}
\tag{Basic}\label{Basic}t\prec q\vee s\ \text{ and }\ t\prec p\vee s\quad&\Rightarrow\quad\exists r\leq p,q\ (t\prec r\vee s).\\
\tag{Basic$'$}\label{Basic'}p\prec q\vee s\ \text{ and }\ q\prec p\vee s\quad&\Rightarrow\quad\exists r\leq p,q\ (p\vee q\leq r\vee s).
\end{align}
\end{prp}

\begin{proof}
By \eqref{pqr}, $p\prec q\vee s$ and $q\prec p\vee s$ implies $p\vee q\prec q\vee s,p\vee s$.  Then \eqref{Basic'} implies $p\vee q\leq r\vee s$, showing that \eqref{Basic} holds.

Conversely, say $S$ is round and \eqref{Basic'} holds.  If $t\prec p\vee s$ then we have finite $E\subseteq S$ such that $E\cup\{t\}$ is far and $q\prec e\vee p\vee s$, for all $e\in E$.  Likewise, $t\prec q\vee s$ yields finite $F\subseteq S$ such that $F\cup\{t\}$ is far and $p\prec f\vee q\vee s$, for all $f\in F$.  Thus, for any $e\in E$ and $f\in F$, we have $p\prec q\vee s\vee e\vee f$ and $q\prec p\vee s\vee e\vee f$.  Then \eqref{Basic'} yields $r_{e,f}\leq p,q$ with $p\vee q\leq r_{e,f}\vee s\vee e\vee f$.  Taking $r=\bigvee_{e\in E, f\in F}r_{e,f}\leq p,q$, we see that $t\prec p\vee q\vee s\leq r\vee s\vee e\vee f$, for all $e\in E$ and $f\in F$.  As $E\cup\{t\}$ is far, \autoref{FarBelow} then yields $t\prec r\vee s\vee f$, for all $f\in F$.  As $F\cup\{t\}$ is far, \autoref{FarBelow} again yields $t\prec r\vee s$, as required.
\end{proof}

\begin{prp}\label{Basic=>Basis}
If $S$ is round and basic then $(\widehat{S}_p)_{p\in S}$ is a basis for the spectrum.
\end{prp}

\begin{proof}
It suffices to show that each $G\in\widehat{S}$ is a filter.  Taking any $e,f,g\in G$, \autoref{MinimalGrillEquivalent} and the fact that $S$ is round yields $e',f'\in S\setminus G$ with $g\prec e\vee e',f\vee f'$ and hence $g\prec e\vee e'\vee f',f\vee e'\vee f'$.  Then \eqref{Basic} yields $d\leq e,f$ with $g\prec d\vee e'\vee f'$.  As $G$ is a grill, $e'\vee f'\notin G$ and hence $d\in G$, showing that $G$ is a filter.
\end{proof}

We can also replace the finite $F$ in the definition of $\prec$ by a singleton.

\begin{thm}
If $S$ is round and basic then, for all $p,q\in S$,
\begin{equation}\label{Rather1}
p\prec q\qquad\Leftrightarrow\qquad\forall s\nleq q\ \exists f\in S\ (\{f,p\}\text{ is far and }s\leq f\vee q).
\end{equation}
\end{thm}

\begin{proof}
$\Leftarrow$ is immediate from the definition of $\prec$.  Conversely, say $p\prec q$.  As $S$ is $\succ$-round, for any $s\nleq q$, we have $t\succ s$.  Then $t\nleq q$ so we have finite $F\subseteq S$ such that $F\cup\{p\}$ is far and $s\prec t\leq f\vee q$, for all $f\in F$.  Successive applications of \eqref{Basic} then yield $e\in S$ with $s\prec e\vee q$ and $e\leq f$, for all $f\in F$, so $\{e,p\}$ is far.
\end{proof}

Lastly we note that \eqref{Basic} implies a version of distributivity for $\prec$.

\begin{prp}\label{precDistributive}
If $S$ is basic then, for all $p,p',s,t\in S$,
\[\tag{$\prec$-Distributive}p'\prec p\leq s\vee t\qquad\Rightarrow\qquad\exists s'\leq s\ \exists t'\leq t\ (p'\prec s'\vee t'\leq p).\]
\end{prp}

\begin{proof}
If $p'\prec p\leq s\vee t$ then $p'\prec p\vee t$, so \eqref{Basic} yields $s'\leq s,p$ with $p'\prec s'\vee t$.  Also $p'\prec s'\vee p$ so \eqref{Basic} again yields $t'\leq t,p$ with $p'\prec s'\vee t'\leq p$.
\end{proof}

\section{Subfitness}\label{Subfitness}

\begin{dfn}
We call $S$ \emph{subfit} if
\[p\nleq q\qquad\Leftrightarrow\qquad\forall p'\geq p\ \exists q'\geq q\ (p\vee q'\geq p'\nleq q').\]
\end{dfn}

\begin{rmk}
If $S$ has a maximum $1$, it again suffices to consider $p'=1$ above, i.e.
\[p\nleq q\qquad\Leftrightarrow\qquad\exists q'\geq q\ (p\vee q'\geq1\nleq q'),\]
which then agrees with subfitness defined for frames in \cite[V.1.1]{PicadoPultr2012}.  The order dual of this was originally called the disjunction property in \cite[Lemma 3]{Wallman1938} and has also been given various other names, e.g. `section semicomplemented' in lattice theory (see \cite{MaedaMaeda1970}), or `separative' in set theory (see \cite{Kunen1980}).
\end{rmk}

\begin{prp}\label{T1subfit}
Any relatively compact $\cup$-subbasis $S$ of a $T_1$ space $X$ is subfit.
\end{prp}

\begin{proof}
If $p\nsubseteq q$ then we have $x\in p\setminus q$.  As $S$ is a subbasis and $X$ is $T_1$, for every $y\in X\setminus p$, we have $r\in S$ with $y\in r\not\ni x$.  For any $p'\supseteq p$, the relative compactness of $p'$ yields finitely many such $r$ covering $p'\setminus p$.  As $S$ is $\cup$-closed we can take the union to obtain a single such $r$ and let $q'=q\cup r\not\ni x$ so $p\cup q'\supseteq p'\nsubseteq q'$, as $x\in p'\setminus q'$.
\end{proof}

\begin{prp}\label{Faithful}
If $S$ is subfit and round then, for any $p,q\in S$,
\[p\leq q\qquad\Leftrightarrow\qquad\widehat{S}_p\subseteq\widehat{S}_q.\]
\end{prp}

\begin{proof}
The $\Rightarrow$ part is immediate.  Conversely, say $p\nleq q$.  As $S$ is round, we have $p'\succ p$.  As $S$ is subfit, we have $q'\geq q$ such that $p\vee q'\geq p'\nleq q'$.  Then $S\setminus q'^\geq$ is a grill containing $p$ but avoiding $q'$.  By \autoref{pGrill}, we have a grill $H\subseteq S\setminus q'^\geq$ containing $p$ with only near finite subsets.  We claim that any non-empty grill $G\subseteq H$ must contain $p'$.  This is immediate if $p'$ is a maximum of $S$.  Otherwise, taking any $g\in G$, the definition of $\prec$ yields finite $E$ such that $E\cup\{p\}$ is far and $g\leq e\vee p'$, for all $e\in E$.  As $H$ only contains near finite subsets, $e\notin H\supseteq G$, for some $e\in E$, and hence $p'\in G$, as $G$ is a grill, proving the claim.  Now Kuratowski-Zorn yields a minimal non-empty grill $G\subseteq H$, necessarily containing $p'$.  Thus $p\in G$ too, as $p'\leq p\vee q'$ and $G$ is a grill containing $p'$ and avoiding $q'$, i.e. $G\in\widehat{S}_p\setminus\widehat{S}_q$.
\end{proof}

In subfit round $\vee$-semilattices, near subsets can be characterised via the spectrum in an analogous manner to \autoref{SubbasisNear} and \autoref{DistributiveNear}.

\begin{prp}\label{NearIntersection}
If $S$ is subfit and round then, for any finite $F\subseteq S$,
\[F\text{ is near}\qquad\Leftrightarrow\qquad\bigcap_{f\in F}\widehat{S}_f\neq\emptyset.\]
\end{prp}

\begin{proof}
The $\Leftarrow$ part was already proved in \autoref{NearMinimalGrill}.  Conversely, say $F$ is near, so we have some $p,q\in S$ with $q\ngeq p\leq q\vee f$, for all $f\in F$.  By \autoref{Faithful}, we have $G\in\widehat{S}_p\setminus\widehat{S}_q$, i.e. $p\in G\not\ni q$ so $F\subseteq G$, as $G$ is a grill, i.e. $G\in\bigcap_{f\in F}\widehat{S}_f$.
\end{proof}

We can then use this to characterise closures in the spectrum.

\begin{cor}\label{Closures}
If $S$ is subfit and round then, for any $p\in S$,
\begin{equation}\label{ClosureChar}
\mathrm{cl}(\widehat{S}_p)=\{G\in\widehat{S}:\forall\text{ finite }F\subseteq G\ (F\cup\{p\}\text{ is near})\}.
\end{equation}
\end{cor}

\begin{proof}
If $G\in\mathrm{cl}(\widehat{S}_p)$ then, as any finite $F\subseteq G$ determines a neighbourhood of $G$, we must have some $H\in\widehat{S}_p\cap\bigcap_{f\in F}\widehat{S}_f$.  But this means $F\cup\{p\}\subseteq H$ and hence $F\cup\{p\}$ is near, by \autoref{NearMinimalGrill}.  Conversely, if $F\cup\{p\}$ is near then $\widehat{S}_p\cap\bigcap_{f\in F}\widehat{S}_f\neq\emptyset$, by \autoref{NearIntersection}, so if this holds for all finite $F\subseteq G$ then $G\in\mathrm{cl}(\widehat{S}_p)$.
\end{proof}

This yields a characterisation of $\prec$ like in \autoref{Faithful}.

\begin{prp}\label{precCharProp}
If $S$ is subfit and round then
\begin{equation}\label{precChar}
p\prec q\qquad\Leftrightarrow\qquad\mathrm{cl}(\widehat{S}_p)\subseteq\widehat{S}_q.
\end{equation}
\end{prp}

\begin{proof}
If $p\nprec q$ then we have $G\in\mathrm{cl}(\widehat{S}_p)\setminus\widehat{S}_q$ by \autoref{pqGrill} and \eqref{ClosureChar}.

Conversely, say we have $G\in\mathrm{cl}(\widehat{S}_p)\setminus\widehat{S}_q$.  In particular, $q$ can not be a maximum of $S$ so, for any $g\in G$, $p\prec q$ would yield finite $F\subseteq S$ such that $F\cup\{p\}$ is far and $g\leq f\vee q$, for all $f\in F$.  As $q\notin G$ and $G$ is a grill, this yields $F\subseteq G$.  As $G\in\mathrm{cl}(\widehat{S}_p)$, $F\cup\{p\}$ is near, by \autoref{Closures}, a contradiction.  Thus $p\nprec q$.
\end{proof}

And now we can finally prove that $\widehat{S}$ is locally closed compact.

\begin{thm}\label{CompactClosures}
If $S$ is subfit and round then $\mathrm{cl}(S_p)$ is compact, for all $p\in S$.
\end{thm}

\begin{proof}
By the Alexander-Wallman subbasis lemma, it suffices to show that every subbasic cover of $\mathrm{cl}(S_p)$ has a finite subcover.  Equivalently, given any ideal $I\subseteq S$ such that $\mathrm{cl}(\widehat{S}_p)\nsubseteq\widehat{S}_j$, for all $j\in I$, we must show that $\mathrm{cl}(\widehat{S}_p)\nsubseteq\bigcup_{j\in I}\widehat{S}_j$.  To see this, first note that $p\nprec j$, for all $j\in I$, by \eqref{precChar}.  As in the proof of \autoref{pqGrill}, we then obtain $G\in\widehat{S}$ such that $G\cap I=\emptyset$ and $F\cup\{p\}$ is near, for all finite $F\subseteq G$.  By \eqref{ClosureChar}, this means $G\in\mathrm{cl}(\widehat{S}_p)\setminus\bigcup_{j\in I}\widehat{S}_j$, as required.
\end{proof}

Our results can be summarised as a duality of the following classes.
\begin{align*}
\mathbf{{}_\cup Sub}&=\text{relatively compact $T_1$ $\cup$-subbases of locally closed compact spaces}.\\
\mathbf{{}_\vee Semi}&=\text{subfit round $\vee$-semilattices}.
\end{align*}

\begin{thm}\label{TheDuality}
$\mathbf{{}_\cup Sub}$ is dual to $\mathbf{{}_\vee Semi}$.
\end{thm}

More precisely, $\mathbf{{}_\cup Sub}\subseteq\mathbf{{}_\vee Semi}$ and the spectrum of any $S\in\mathbf{{}_\cup Sub}$ recovers the original space in which $S$ lies, by \autoref{Recovery}.  Conversely every $S\in\mathbf{{}_\vee Semi}$ has a locally closed compact $T_1$ spectrum on which $S$ is faithfully represented as $(\widehat{S}_p)_{p\in S}\in\mathbf{{}_\cup Sub}$, by \autoref{MinimalT1}, \autoref{Faithful} and \autoref{CompactClosures}.

By \autoref{Basic=>Basis}, we could restrict \autoref{TheDuality} to obtain a duality between
\begin{align*}
\mathbf{{}_\cup Basis}&=\text{relatively compact $T_1$ $\cup$-bases of locally closed compact spaces}.\\
\mathbf{{}_\vee BSemi}&=\text{basic subfit round $\vee$-semilattices}.
\end{align*}
Moreover, $\prec$ has a truly first order equivalent definition in $\mathbf{{}_\vee BSemi}$ (see \eqref{Rather1}) so $\mathbf{{}_\vee BSemi}$ forms an elementary class in the usual model theoretic sense.

We could also restrict \autoref{TheDuality} to
\begin{align*}
\mathbf{{}_\cup SubC}&=\text{$T_1$ $\cup$-bases of compact spaces}.\\
\mathbf{{}_\vee Semi1}&=\text{unital subfit $\vee$-semilattices}.
\end{align*}
As with Wallman's work, our primary motivation was to obtain order theoretic duals of topological spaces.  However, this shows that our duality could also be used to obtain topological representations of general unital (i.e. having a maximum $1$) subfit $\vee$-semilattices.  Topological representations of more general $\vee$-semilattices have been obtained in \cite{JipsenMoshier2014}, but only on highly non-$T_1$ `HMS-spaces'.

\section{Functoriality}\label{Functoriality}

We now set about making our duality functorial.  Our morphisms in $\mathbf{{}_\vee Semi}$ will be relations, while our morphisms in $\mathbf{{}_\cup Sub}$ will be (partial) functions.  As usual, we define the composition $\phi'\circ\phi$ of relations $\phi\subseteq X'\times X$ and $\phi'\subseteq X''\times X'$ by
\[\phi'\circ\phi=\{(x'',x)\in X''\times X:\exists x'\in X'\ ((x'',x')\in\phi'\text{ and }(x',x)\in\phi\}.\]
We denote the image of any $Y\subseteq X$ and preimage of any $Y'\subseteq X'$ under $\phi$ by
\begin{align*}
\phi[Y]&=\{x'\in X':\exists y\in Y\ ((x',y)\in\phi)\}.\\
[Y']\phi&=\{x\in X:\exists y'\in Y'\ ((y',x)\in\phi)\}.
\end{align*}
If $\phi[\{x\}]$ contains at most one element, for each $x\in X$, then we view $\phi$ as a function from a subset of $X$ to a subset of $X'$ with the usual domain and range denoted by
\[\mathrm{dom}(\phi)=[X']\phi\qquad\text{and}\qquad\mathrm{ran}(\phi)=\phi[X].\]

Given a subbasis $S$ of $X$, let us call unions of elements of $S$ \emph{wide open}, i.e.
\[O\subseteq X\emph{ is wide open}\qquad\Leftrightarrow\qquad\exists Q\subseteq S\ (O=\bigcup Q).\]

\begin{dfn}
Given subbases $S$ and $S'$ of $X$ and $X'$, a function $\phi\subseteq X\times X'$ is
\begin{enumerate}
\item \emph{wide continuous} if $[O']\phi$ is wide open whenever $O'\subseteq X'$ is.
\item \emph{closed compact} if $\phi[C]$ is closed compact whenever $C\subseteq\mathrm{dom}(\phi)$ is.
\end{enumerate}
We call $\phi$ \emph{very continuous} if $\phi$ is both wide continuous and closed compact.
\end{dfn}

In other words, a function is wide continuous if preimages of wide open sets are wide open, while a map is closed compact if (forward) images of closed compact subsets of the domain are closed compact.  We emphasise that we are allowing partial functions here, i.e. $\mathrm{dom}(\phi)$ does not have to be the entirety of $X$ (although it does have to be wide open if $\phi$ is wide continuous and $S'$ covers $X'$).

If $S$ is a basis of $X$ then any continuous $\phi\subseteq X'\times X$ is automatically wide continuous.  If $X'$ is Hausdorff then any continuous $\phi\subseteq X'\times X$ is automatically closed compact.  So if we restrict to bases of Hausdorff spaces then very continuous maps are just the usual continuous maps.  In general, however, this is not so, e.g. if $X=[0,1]$ with the subbasis $S=\{(r,1]:r\in\mathbb{R}\}\cup\{[0,r):r\in\mathbb{R}\}$ then, for any function $\phi$ on $X$,
\[\phi\text{ is very continuous}\qquad\Leftrightarrow\qquad\phi\text{ is monotone continuous}\]
(where \emph{monotone} means betweenness preserving, i.e. order preserving or reversing).

\begin{prp}
$\mathbf{{}_\cup Sub}$ forms a category with very continuous morphisms.
\end{prp}

\begin{proof}
Just note that very continuous maps are closed under composition and also include identity maps, which are immediately seen to be identity morphisms.
\end{proof}

To describe the morphisms in $\mathbf{{}_\vee Semi}$, it will be convenient to introduce the formal expression $\bigwedge E\leq\bigwedge F$ to mean that, for all $p,q\in S$,
\[\forall e\in E\ (p\leq q\vee e)\quad\Rightarrow\quad\forall f\in F\ (p\leq q\vee f).\]

\begin{prp}
For all finite $E,F\subseteq S$,
\begin{equation}\label{bigwedges}
\bigwedge E\leq\bigwedge F\qquad\Rightarrow\qquad\bigcap_{e\in E}\widehat{S}_e\subseteq\bigcap_{f\in F}\widehat{S}_f.
\end{equation}
The converse also holds if $S$ is round and subfit.
\end{prp}

\begin{proof}
Say $\bigwedge E\leq\bigwedge F$ and take any $G\in\bigcap_{e\in E}\widehat{S}_e$.  For any $p\in G$ and $e\in E\subseteq G$, \autoref{MinimalGrillEquivalent} yields $q_e\in S\setminus G$ with $p\leq q_e\vee e$.  Taking $q=\bigvee_{e\in E}q_e\notin G$, as $G$ is a grill, we see that $p\leq q\vee e$, for all $e\in E$.  For all $f\in F$, it follows that $p\leq q\vee f$ and hence $f\in G$, as $\bigwedge E\leq\bigwedge F$ and $G$ is a grill, i.e. $G\in\bigcap_{f\in F}\widehat{S}_f$.

Conversely, say $S$ is round and subfit and $\bigcap_{e\in E}\widehat{S}_e\subseteq\bigcap_{f\in F}\widehat{S}_f$.  If $p\leq q\vee e$, for all $e\in E$, then $\widehat{S}_p\subseteq\widehat{S}_q\cup\widehat{S}_e$, for all $e\in E$, and hence $\widehat{S}_p\subseteq\widehat{S}_q\cup\bigcap_{e\in E}\widehat{S}_e\subseteq\widehat{S}_q\cup\bigcap_{f\in F}\widehat{S}_e$, i.e. $\widehat{S}_p\subseteq\widehat{S}_q\cup\widehat{S}_f=\widehat{S}_{q\vee f}$ and hence $p\leq q\vee f$, for all $f\in F$, by \eqref{Faithful}.
\end{proof}

Given a relation $\sqsubset\ \subseteq S\times S'$ and $T\subseteq S$, let
\[T^\sqsubset=[T]\!\sqsubset\ =\{t'\in S':\exists t\in T\ (t\sqsubset t')\}.\]
Extend $\sqsubset$ and its opposite $\sqsupset$ to subsets $T\subseteq S$ and $T'\subseteq S'$ by defining
\begin{align*}
T\sqsubset T'\qquad&\Leftrightarrow\qquad\forall t\in T\ \exists t'\in T'\ (t\sqsubset t')\qquad\Leftrightarrow\qquad T\subseteq T'^\sqsupset\\
T\sqsupset T'\qquad&\Leftrightarrow\qquad\forall t\in T\ \exists t'\in T'\ (t\sqsupset t')\qquad\Leftrightarrow\qquad T\subseteq T'^\sqsubset.
\end{align*}

\begin{dfn}
For any $S,S'\in{}_\vee\mathbf{Semi}$, we call $\sqsubset\ \subseteq S\times S'$ a \emph{$\vee$-relation} if
\begin{enumerate}
\myitem[(Auxiliary)]\label{Auxiliary} For all $p,q\in S$ and $p',q'\in S'$,
\[p\leq q\sqsubset q'\leq p'\quad\Rightarrow\quad p\sqsubset p'.\]
\myitem[($\vee$-Preserving)]\label{veePreserving} For all $p,q\in S$ and $r'\in S'$,
\[p,q\sqsubset r'\quad\Rightarrow\quad p\vee q\sqsubset r'.\]
\myitem[(Decomposition)]\label{Decomposition} For all $p,q\in S$ and finite $F'\subseteq S'$,
\[p\vee q\succ p\sqsubset\bigvee F'\quad\Rightarrow\quad\exists\text{ finite }F\sqsubset F'\ (p\leq q\vee\bigvee F).\]
\myitem[(Complementation)]\label{Complementation} For all $p,q\in S$ and $p',q'\in S'$ with $p\vee q\succ p\sqsubset p'$, we have finite $F'\subseteq S'$ with $q'\leq p'\vee f'$, for all $f'\in F'$, and
\[F'\sqsupset F\quad\Rightarrow\quad\bigwedge(F\cup\{p\})\leq q.\]
\end{enumerate}
\end{dfn}

See \eqref{subphi} and \autoref{cts->vmorph} below for what motivates these conditions.

\begin{rmk}
To motivate \ref{Complementation} in particular, imagine that $S$ and $S'$ are concrete subbases of spaces $X$ and $X'$ and that $\sqsubset\ =\ \sqsubset_\phi$ is defined from a very continuous function $\phi$ as in \eqref{subphi} below.  What \ref{Complementation} is essentially saying then is that the image under $\phi$ of the closed set $p\setminus q$ is contained in the closed set $p'\setminus\bigcap F'$, which corresponds to the fact $\phi$ is compact closed.
\end{rmk}

\begin{thm}\label{veeRelations->Maps}
If $\sqsubset$ is a $\vee$-relation, the map $\phi_\sqsubset\subseteq\widehat{S'}\times\widehat{S}$ is very continuous, where
\[\mathrm{dom}(\phi_\sqsubset)=\{G\in\widehat{S}:G^\sqsubset\neq\emptyset\}\quad\text{and}\quad\phi_\sqsubset(G)=G^\sqsubset.\]
\end{thm}

\begin{proof}
First we show that $\mathrm{ran}(\phi_\sqsubset)\subseteq\widehat{S'}$, i.e. $G^\sqsubset\in\widehat{S}'$ whenever $G\in\widehat{S}$ and $G^\sqsubset\neq\emptyset$.

By \ref{Auxiliary}, $G^\sqsubset$ is an up-set.  To see that $G^\sqsubset$ is a grill, say $G\ni p\sqsubset r'\vee s'$.  As $S$ is round, we have $t\succ p$ and then \autoref{MinimalGrillEquivalent} yields $q\in S\setminus G$ with $p\prec t\leq p\vee q$.  Then \ref{Decomposition} yields finite $F\sqsubset\{r',s'\}$ with $p\leq q\vee\bigvee F$.  As $q\notin G\ni p$ and $G$ is a grill, $F\cap G\neq\emptyset$ and hence $\{r',s'\}\cap G^\sqsubset\neq\emptyset$, showing that $G^\sqsubset$ is also a grill.

To see that $G^\sqsubset$ is proper and minimal, say $G\ni p\sqsubset p'$ and $q'\in S'$.  Again we have $q\in S\setminus G$ with $p\prec p\vee q$.  Then \ref{Complementation} yields finite $F'\subseteq S'$ with $q'\leq p'\vee f'$, for all $f'\in F'$, and $\bigwedge(F\cup\{p\})\leq q$ whenever $F'\sqsupset F$.  If we had $F'\subseteq G^\sqsubset$ then we would indeed have finite $F\subseteq G$ with $F'\sqsupset F$ and hence $\bigwedge(F\cup\{p\})\leq q$.  But then \eqref{bigwedges} yields $G\in\widehat{S}_p\cap\bigcap_{f\in F}\widehat{S}_f\subseteq\widehat{S}_q$, i.e. $q\in G$, a contradiction.  Thus $F'\nsubseteq G^\sqsubset$, i.e. we have $f'\in F'\setminus G^\sqsubset$ with $q'\leq p'\vee f'$ so $G'$ is proper and minimal, by \autoref{MinimalGrillEquivalent}.  This shows that $\mathrm{ran}(\phi_\sqsubset)\subseteq\widehat{S'}$.

To see that $\phi_\sqsubset$ is wide continuous, just note that, for all $p'\in S'$,
\[[\widehat{S}'_{p'}]\phi_\sqsubset=\bigcup\limits_{p\sqsubset p'}\widehat{S}_p.\]
To see that $\phi_\sqsubset$ is closed compact, take any closed compact $C\subseteq\mathrm{dom}(\phi_\sqsubset)$.  As $\phi_\sqsubset$ is (wide) continuous, we know that $\phi_\sqsubset[C]$ is compact.  Say $\phi_\sqsubset[C]$ were not closed, so we have some $G'\in\mathrm{cl}(\phi_\sqsubset[C])\setminus\phi_\sqsubset[C]$.  For every $G\in C$, this means that $G'\neq G^\sqsubset$ so we have some $p'\in G^\sqsubset\setminus G$, and thus we have some $p\in G$ with $p\sqsubset p'$.  As $C$ is compact, finitely many such $p$ cover $C$ so, taking joins and using \ref{veePreserving}, we get $p\in G$ with $p\sqsubset p'\notin G$ and $C\subseteq\widehat{S}_p$.  As $C$ is a closed subset of $\widehat{S}_p$, we have finite $E\subseteq S$ with $p\prec e\vee p$, for all $e\in E$, and $\bigcap_{e\in E}\widehat{S}_e\cap C=\emptyset$ (each point in $\mathrm{cl}(p)\setminus p$ has a basic neighbourhood $\bigcap_{e\in D}\widehat{S}_d$ disjoint from $C$, so we can cover $C$ by finitely many such neighbourhoods defined by finite $E_1,\cdots,E_n$ and let $E=E_1\vee\cdots\vee E_n$).

Taking any $g'\in G'$, for each $e\in E$, \ref{Complementation} yields finite $F'_e\subseteq S'$ with $g'\leq p'\vee f'$, for all $f'\in F'_e$, and $\bigwedge(F\cup\{p\})\leq e$ whenever $F'_e\sqsupset F$.  Taking $F'=\bigcup_{e\in E}F'_e$, it follows that $g'\leq p'\vee f'$, for all $f'\in F'$, and $\bigwedge(F\cup\{p\})\leq\bigwedge E$ whenever $F'\sqsupset F$.  As $p'\notin G'$, $g'\in G'$ and $G'$ is a grill, it follows that $F'\subseteq G'$ so $G'\in\bigcap_{f'\in F'}\widehat{S}_{f'}$.  Thus $\phi_\sqsubset[C]\cap\bigcap_{f'\in F'}\widehat{S}_{f'}\neq\emptyset$, as $G'\in\mathrm{cl}(\phi_\sqsubset[C])$.  This means we have $H\in C$ with $F'\subseteq H^\sqsubset$ and hence we have finite $F\subseteq H$ with $F'\sqsupset F$.  Then $\bigwedge(F\cup\{p\})\leq\bigwedge E$ so, as $H\in C\subseteq\widehat{S}_p$ and hence $p\in H$ too, \eqref{bigwedges} yields $H\in\bigcap_{e\in E}\widehat{S}_e\cap C=\emptyset$, a contradiction.  This shows that $\phi_\sqsubset$ is closed compact and hence very continuous.
\end{proof}

\begin{thm}\label{veeSemiCat}
${}_\vee\mathbf{Semi}$ forms a category with $\vee$-relations as morphisms.
\end{thm}

\begin{proof}
Note that $\leq$ itself is always a $\vee$-relation on any $S\in{}_\vee\mathbf{Semi}$.  Moreover, \ref{Auxiliary} is saying that $\sqsubset$ coincides with $\leq\circ\sqsubset$ and $\sqsubset\circ\leq$, i.e. $\leq$ is always an identity among $\vee$-relations.  So we just have to show $\vee$-relations are closed under composition, i.e. if $\sqsubset\ \subseteq S\times S'$ and $\sqsubset'\ \subseteq S'\times S''$ are $\vee$-relations then so is $\sqsubset\circ\sqsubset'$.

To see this, first note that if $s\leq t\sqsubset\circ\sqsubset't''\leq s''$, i.e. $s\leq t\sqsubset t'\sqsubset't''\leq s''$, for some $s'\in S'$, then $s\sqsubset t'\sqsubset's''$, by \ref{Auxiliary}, and hence $s\sqsubset\circ\sqsubset's''$, by the definition of relation composition.  This shows that $\sqsubset\circ\sqsubset'$ satisfies \ref{Auxiliary}.  We also immediately see that $\sqsubset\circ\sqsubset'$ satisfies \ref{veePreserving}, as $\sqsubset$ and $\sqsubset'$ do.

To show that $\sqsubset\circ\sqsubset'$ satisfies \ref{Complementation}, say $p\vee q\succ p\sqsubset p'\sqsubset'p''$ and $q''\in S''$.  Take $q'\in S'$ with $p'\prec q'$ and take finite $F'\in S'$ witnessing \ref{Complementation} for $\sqsubset$ so $p'\prec q'\leq p'\vee f'$, for all $f'\in F'$, and $\bigwedge(F\cup\{p\})\leq q$ whenever $F'\sqsupset F$.  Now take finite $E''\subseteq S''$ witnessing \ref{Complementation} for $\sqsubset'$, for all $f'\in F'$ like in the above proof, so $q''\leq p''\vee e''$, for all $e''\in E''$, and $\bigwedge(E'\cup\{p'\})\leq\bigwedge F'$ whenever $E''\sqsupset'E'$.  We claim that $E''$ also witnesses \ref{Complementation} for $\sqsubset\circ\sqsubset'$, i.e. $\bigwedge(E\cup\{p\})\leq q$ whenever $E''\sqsupset'E'\sqsupset E$.

To see this, say we have finite $E'$ and $E$ with $E''\sqsupset E'\sqsupset E$ but $\bigwedge(E\cup\{p\})\nleq q$.  Then we have $G\in\bigcap_{e\in E}\widehat{S}_e\cap\widehat{S}_p\setminus\widehat{S}_q$, by \eqref{bigwedges}, i.e. $E\cup\{p\}\subseteq G\not\ni q$.  It follows that $E'\cup\{p'\}\subseteq G^\sqsubset$ so \eqref{bigwedges} again yields $F'\subseteq G^\sqsubset$, as $\bigwedge(E'\cup\{p'\})\leq\bigwedge F'$ and $G^\sqsubset\in\widehat{S'}$, by \autoref{FunctorTheorem}.  But then we would have finite $F\subseteq G$ with $F'\sqsupset F$ so \eqref{bigwedges} yet again yields $q\in G$, as $\bigwedge(F\cup\{p\})\leq q$, a contradiction.

To show that $\sqsubset\circ\sqsubset'$ satisfies \ref{Decomposition}, assume above that $p''=\bigvee D''$, for some finite $D''$.  For each $f'\in F'$, \ref{Decomposition} for $\sqsubset'$ yields finite $D'_{f'}\sqsubset'D''$ with $p'\leq f'\vee\bigvee D'_{f'}$ and then \ref{Decomposition} for $\sqsubset$ yields finite $D_{f'}\sqsubset\{f'\}\cup D'_{f'}$ with $p\leq q\vee\bigvee D_{f'}$.  Letting $D=\bigcup_{f'\in F'}(D_{f'}\cap D_{f'}'^\sqsupset)$, so $D\sqsubset\circ\sqsubset'D''$, we claim $p\leq q\vee\bigvee D$.  If not, we would have $G\in\widehat{S}_p\setminus(\widehat{S}_q\cup\bigcup_{d\in D}\widehat{S}_d)$.  For each $f'\in F'$, $p\leq q\vee\bigvee D_{f'}$ so we have some $g_{f'}\in G\cap D_{f'}\subseteq D_{f'}\setminus D\subseteq D_{f'}\setminus D_{f'}'^\sqsupset$ and hence $g_{f'}\sqsubset f'$, as $D_{f'}\sqsubset\{f'\}\cup D'_{f'}$.  Taking $F=\{g_{f'}:f'\in F'\}\subseteq G$, we see that $F'\sqsupset F$ so $\bigwedge(F\cup\{p\})\leq q$ and hence $q\in G$, by \eqref{bigwedges}, a contradiction.  This proves the claim, which shows that $\sqsubset\circ\sqsubset'$ satisfies \ref{Decomposition} and is thus a $\vee$-relation.
\end{proof}

\begin{cor}\label{FunctorTheorem}
We have a functor from ${}_\vee\mathbf{Semi}$ to ${}_\cup\mathbf{Sub}$ given by
\[S\mapsto(\widehat{S}_p)_{p\in S}\qquad\text{and}\qquad\sqsubset\ \mapsto\phi_\sqsubset.\]
\end{cor}

\begin{proof}
By \autoref{veeRelations->Maps}, $\sqsubset\ \mapsto\phi_\sqsubset$ takes morphisms to morphisms.  It also preserves composition, as $\phi_{\sqsubset\circ\sqsubset'}(G)=G^{\sqsubset\circ\sqsubset'}=(G^\sqsubset)^{\sqsubset'}=\phi_{\sqsubset'}(\phi_\sqsubset(G))=(\phi_{\sqsubset'}\circ\phi_\sqsubset)(G)$.\qedhere
\end{proof}

Given subbases $S$ and $S'$ on spaces $X$ and $X'$ respectively and a function $\phi$ with $\mathrm{dom}(\phi)\subseteq X$ and $\mathrm{ran}(\phi)\subseteq X'$, we define a relation $\sqsubset_\phi\ \subseteq S\times S'$ by
\begin{equation}\label{subphi}
s\sqsubset_\phi s'\qquad\Leftrightarrow\qquad s\subseteq[s']\phi.
\end{equation}

\begin{thm}\label{cts->vmorph}
For any $S,S'\in{}_\cup\mathbf{Sub}$ on $X$ and $X'$ and any very continuous $\phi\subseteq X'\times X$, the relation $\sqsubset_\phi\ \subseteq S\times S'$ is a $\vee$-relation and, for all $x\in\mathrm{dom}(\phi)$,
\begin{equation}\label{p=psp}
S'_{\phi(x)}=\phi_{\sqsubset_\phi}(S_x).
\end{equation}
\end{thm}

\begin{proof}
If $p\subseteq q\subseteq\phi^{-1}[q']$ and $q'\subseteq p'$ then $p\subseteq\phi^{-1}[q']\subseteq\phi^{-1}[p']$, showing that $\sqsubset_\phi$ satisfies \ref{Auxiliary}.  If $p\subseteq\phi^{-1}[p']$ and $q\subseteq\phi^{-1}[q']$ then $p\cup q\subseteq\phi^{-1}[p']\cup\phi^{-1}[q']=\phi^{-1}[p'\cup q']$, showing that $\sqsubset_\phi$ also satisfies \ref{veePreserving}.

Now say $q'\in S'$ and $p\vee q\succ p\subseteq[\bigcup F']\phi$, for some finite $F'\subseteq S'$.  By \autoref{T1round}, $\mathrm{cl}(p)\subseteq p\cup q$ so $p\setminus q=\mathrm{cl}(p)\setminus q\subseteq\mathrm{cl}(p)$ is closed and compact.  As $\phi$ is wide continuous, we can cover $p\setminus q\subseteq p\subseteq\bigcup_{f'\in F'}[f']\phi$ with finitely many $f\in S$ each contained in $[f']\phi$, for some $f'\in F'$.  This yields finite $F\sqsubset_\phi F'$ with $p\subseteq q\cup\bigcup F$, showing that $\sqsubset_\phi$ satisfies \ref{Decomposition}.

Again if $q'\in S'$ and $p\vee q\succ p\subseteq[p']\phi$ then $\phi[p\setminus q]$ is closed, as $\phi$ is closed compact.  As $\phi[p\setminus q]\subseteq\phi[p]\subseteq p'$ is also disjoint from the compact set $\mathrm{cl}(q')\setminus p'$, we can cover the latter with finitely many finite intersections $\bigcap E'_1,\cdots,\bigcap E'_n$ from $S'$ disjoint from $\phi[p\setminus q]$.  As $S'$ is a $\cup$-basis, we have finite $F'\subseteq S'$ defined by
\[F'=\{\bigcup E':E'\subseteq\bigcup_kE'_k\text{ and }\forall k\ (E'\cap E'_k\neq\emptyset)\}.\]
Then $\bigcap F'=\bigcup_k\bigcap E_k$ is disjoint from $\phi[p\setminus q]$ and covers $\mathrm{cl}(q')\setminus p'$ and hence $q'\subseteq p'\cup f'$, for all $f'\in F'$.  Thus $[\bigcap F']\phi$ is disjoint from $p\setminus q$ and hence the same is true for $\bigcap F$ whenever $F'\sqsupset_\phi F$, as this implies $\bigcap F\subseteq\bigcap_{f'\in F'}[f']\phi\subseteq[\bigcap F']\phi$.  In other words, $p\cap\bigcap F\subseteq q$ and hence $\bigwedge(F\cup\{p\})\leq q$, showing that $\sqsubset_\phi$ satisfies \ref{Complementation}.

Lastly for \eqref{p=psp} note that, for any $x\in\mathrm{dom}(\phi)$, the wide continuity of $\phi$ yields
\[S_x^{\sqsubset_\phi}=\{s'\in S':\exists s\in S\ (x\in s\subseteq[s']\phi)\}=\{s'\in S':\phi(x)\in s'\}=S'_{\phi(x)}.\qedhere\]
\end{proof}

It follows that the functor in \autoref{FunctorTheorem} is full, but it is not faithful, as we always have $\phi_\sqsubset=\phi_\sqsubseteq$ where $\sqsubseteq$ is the weakening of $\sqsubset$ defined by
\[p\sqsubseteq p'\qquad\Leftrightarrow\qquad\forall q\in S\ (p\prec p\vee q\ \Rightarrow\ \exists r\sqsubset p'\ (p\leq q\vee r)).\]

\begin{prp}
If $\sqsubset\ \subseteq S\times S'$ is a $\vee$-relation then
\[p\sqsubseteq p'\qquad\Leftrightarrow\qquad\widehat{S}_p\sqsubset_{\phi_\sqsubset}\widehat{S'}_{p'}.\]
\end{prp}

\begin{proof}
First note that
\[\widehat{S}_p\sqsubset_{\phi_\sqsubset}\widehat{S'}_{p'}\quad\Leftrightarrow\quad\widehat{S}_p\subseteq[S'_{p'}]\phi_\sqsubset\quad\Leftrightarrow\quad\forall G\in\widehat{S}\ (p\in G\ \Rightarrow\ p'\in G^\sqsubset).\]

Now say $G\ni p\sqsubseteq p'$.  As $S$ is round, \autoref{MinimalGrillEquivalent} yields $q\in S\setminus G$ with $p\prec p\vee q$.  As $p\sqsubseteq p'$, we have $r\sqsubset p'$ with $p\leq q\vee r$ and hence $r\in G$, as $G$ is a grill and $p\in G\not\ni q$.  Thus $p'\in G^\sqsubset$, showing that $p\sqsubseteq p'$ implies $\widehat{S}_p\sqsubset_{\phi_\sqsubset}\widehat{S'}_{p'}$.

Conversely, say $\widehat{S}_p\sqsubset_{\phi_\sqsubset}\widehat{S}'_{p'}$ and $p\prec p\vee q$.  By \ref{veePreserving}, $p'^\sqsupset$ is an ideal.  If there were no $r\sqsubset p'$ with $p\leq q\vee r$ then $S\setminus p'^\sqsupset$ would be a grill $G$ containing $p$ but avoiding $q$.  Arguing as in \autoref{Faithful}, we find a minimal grill $H\subseteq G$ with $p\vee q\in H$ and hence $p\in H$, as $q\notin G\supseteq H$.  But then $H$ avoids $p'^\sqsupset$ so $p'\notin H^\sqsubset$, contradicting $\widehat{S}_p\sqsubset_{\phi_\sqsubset}\widehat{S'}_{p'}$.  Thus we must have had finite $r\sqsubset p'$ with $p\leq q\vee r$, showing that $p\sqsubseteq p'$.
\end{proof}

So if we want a categorical equivalence rather than just a functor in \autoref{FunctorTheorem}, these results show that this could be achieved by restricting to $\vee$-relations satsifying $\sqsubset\ =\ \sqsubseteq$ and redefining the composition of such $\sqsubset$ and $\sqsubset'$ as $\underline{\sqsubset\circ\sqsubset'}$.  For a different approach to the locally compact (locally) Hausdorff case, see \cite{BiceStarling2018}.

Also, if one wants to deal with total functions in ${}_\cup\mathbf{Sub}$ then in ${}_\vee\mathbf{Semi}$ one could further require $\vee$-relations to satisfy $S\subseteq S'^\sqsupset$, i.e. $\forall s\in S\ \exists s'\in S'\ (s\sqsubset s')$.

Let us also note that a more algebraic functorialisation of the classic Wallman duality can be found in \cite{BialasBlaszczyk2015} (which is based on earlier work of the second author of the present paper in \cite{Kubis2014}).  Specifically, in \cite{BialasBlaszczyk2015} they showed how to represent lattice homomorphisms (rather than $\vee$-relations) as continuous functions.  Consequently, the resulting functor is far from being full, as it ignores the many continuous maps that do not arise as lattice homomorphisms between fixed bases (and also only applies normal lattices/Hausdorff spaces, not $T_1$ spaces).

\section{Cover Relations}\label{CoverRelations}

Here we briefly consider how one might extend Wallman duality to arbitrary subbases, i.e. without any semilattice structure coming from unions.  The first thing to note is that the inclusion relation on an arbitrary subbasis could easily reduce to mere equality (and thus reveal nothing about the underlying space).  To rectify this we consider `cover relations' on finite subsets of the subbasis.

Denote the finite subsets of a set $S$ by $\mathcal{F}(S)=\{F\subseteq S:F\text{ is finite}\}$.

\begin{dfn}
We call a relation $\looparrowright$ on $\mathcal{F}(S)$ a \emph{cover relation} if
\begin{align}
\label{Antisymmetric+Reflexive}\{p\}\looparrowright\{q\}\quad\text{and}\quad\{q\}\looparrowright\{p\}\qquad&\Leftrightarrow\qquad p=q.\\
\label{Transitive+Monotone}\{p\}\cup D\looparrowright E\quad\text{and}\quad D\looparrowright E\cup\{p\}\qquad&\Leftrightarrow\qquad D\looparrowright E.
\end{align}
\end{dfn}

Alternatively, we could write these properties as four separate axioms as follows.
\begin{align}
\tag{Reflexive}\label{Reflexive}\{p\}\looparrowright\{p\}\qquad&\\
\tag{Antisymmetric}\label{Antisymmetric}\{p\}\looparrowright\{q\}\quad\text{and}\quad\{q\}\looparrowright\{p\}\qquad&\Rightarrow\qquad p=q.\\
\tag{Monotone}\label{Monotone}\{p\}\cup D\looparrowright E\quad\text{and}\quad D\looparrowright E\cup\{p\}\qquad&\Leftarrow\qquad D\looparrowright E.\\
\tag{Transitive}\label{Transitive}\{p\}\cup D\looparrowright E\quad\text{and}\quad D\looparrowright E\cup\{p\}\qquad&\Rightarrow\qquad D\looparrowright E.
\end{align}

The name for \eqref{Transitive} is due to the fact it implies transitivity on singletons, i.e. $\{p\}\looparrowright\{q\}\looparrowright\{r\}$ implies $\{p\}\looparrowright\{r\}$, as we will see in \autoref{Covers->Joins} below.

The canonical example of a cover relation is given is follows.

\begin{prp}
If $S\subseteq\mathcal{O}(X)$ then we can define a cover relation on $S$ by
\[D\looparrowright E\qquad\Leftrightarrow\qquad\bigcap D\subseteq\bigcup E.\]
\end{prp}

\begin{proof}
We check required conditions.
\begin{itemize}
\item[\eqref{Antisymmetric+Reflexive}] Just note $p\subseteq q$ and $q\subseteq p$ iff $p=q$.

\item[\eqref{Transitive+Monotone}] If $p\cap\bigcap D\subseteq\bigcup E$ and $\bigcap D\subseteq\bigcup E\cup p$ then
\[\bigcap D\ \subseteq\ (\bigcup E\cup p)\cap\bigcap D\ \subseteq\ \bigcup E\cup(p\cap\bigcap D)\ \subseteq\ \bigcup E.\]
Conversely, if $\bigcap D\subseteq\bigcup E$ then $p\cap\bigcap D\subseteq\bigcup E$ and $\bigcap D\subseteq\bigcup E\cup p$. \qedhere
\end{itemize}
\end{proof}

Given a cover relation, we would call $G\subseteq S$ a \emph{grill} when
\[\tag{Grill}G\supseteq F\looparrowright E\qquad\Rightarrow\qquad G\cap E\neq\emptyset.\]
Using these grills, we could again obtain a duality with subbases of locally closed compact $T_1$ spaces.  Instead of redoing all the relevant theory, we just show how to pass between $\vee$-semilattices and cover relations.

\begin{prp}\label{Joins->Covers}
If $S$ is a $\vee$-semilattice, a cover relation on $\mathcal{F}(S)$ is given by
\[D\looparrowright E\qquad\Rightarrow\qquad\forall p,q\in S\ (\forall d\in D\ (p\leq q\vee d)\ \Rightarrow\ p\leq q\vee\bigvee E).\]
\end{prp}

\begin{proof}
If $d\leq e$ then $p\leq q\vee d$ implies $p\leq q\vee e$, i.e. $\{d\}\looparrowright\{e\}$.  Conversely, if $\{d\}\looparrowright\{e\}$ then, as $d\leq e\vee d$, we must have $d\leq e\vee e=e$ (taking $p=d$ and $q=e$ in the definition of $\looparrowright$), i.e.
\[d\leq e\qquad\Leftrightarrow\qquad\{d\}\looparrowright\{e\}.\]
In particular, $d=e$ iff $\{d\}\looparrowright\{e\}$ and $\{e\}\looparrowright\{d\}$, i.e. \eqref{Antisymmetric+Reflexive} holds.

Note \eqref{Monotone} is immediate.  Conversely, say $\{r\}\cup D\looparrowright E$ and $D\looparrowright E\cup\{r\}$ and we are given $p,q\in S$ with $p\leq q\vee d$, for all $d\in D$.  Then certainly $p\leq q\vee\bigvee E\vee d$, for all $d\in D$, and $p\leq q\vee\bigvee E\vee r$, as $D\looparrowright E\cup\{r\}$.  Applying the definition of $D\cup\{r\}\looparrowright E$ with $q\vee\bigvee E$ in place of $q$ then yields $p\leq q\vee\bigvee E\vee\bigvee E=q\vee\bigvee E$, showing that $D\looparrowright E$.  Thus \eqref{Transitive} holds and $\looparrowright$ is a cover relation.
\end{proof}

\begin{prp}\label{Covers->Joins}
If $\looparrowright$ is a cover relation, a preorder on $\mathcal{F}(S)$ is given by
\[D\leq E\qquad\Leftrightarrow\qquad\forall d\in D\ (\{d\}\looparrowright E).\]
Identifying $D$ and $E$ whenever $D\leq E\leq D$ then yields a $\vee$-semilattice with $\vee=\cup$.
\end{prp}

\begin{proof}
As \eqref{Reflexive} yields $\{e\}\looparrowright\{e\}$, for all $e\in E$, \eqref{Monotone} then yields $\{e\}\looparrowright E$, for all $e\in E$, and hence $E\leq E$ (and $\emptyset\leq\emptyset$ holds vacuously), i.e. $\leq$ is reflexive.  For transitivity, say $C\leq D\leq E$ so $\{c\}\looparrowright D$, for all $c\in C$, and $\{d\}\looparrowright E$, for all $d\in D$.  Let $D=\{d_1,\ldots,d_n\}$.  By \eqref{Transitive+Monotone}, $c\looparrowright D\cup E$ and $\{c,d_1\}\looparrowright(D\setminus\{d_1\})\cup E$ and hence $c\looparrowright(D\setminus\{d_1\})\cup E$.  Again by \eqref{Transitive+Monotone}, $\{c,d_2\}\looparrowright(D\setminus\{d_1,d_2\})\cup E$ and hence $c\looparrowright(D\setminus\{d_1,d_2\})\cup E$.  Continuing in this way, we obtain $c\looparrowright E$, for all $c\in C$, and hence $C\leq E$.  Thus $\leq$ is indeed a preorder and hence a partial order when we identify $D$ and $E$ whenever $D\leq E\leq D$.  To see that the resulting poset is a $\vee$-semilattice with $\vee=\cup$, just note that $D\leq D\cup E$ and $E\leq D\cup E$, by \eqref{Monotone}, and if $D,E\leq F$ then the definition of $\leq$ immediately yields $D\cup E\leq F$.
\end{proof}

Say we are given a cover relation $\looparrowright$ and we define $\leq$ on $\mathcal{F}(S)$ as above in \autoref{Covers->Joins}.  Let us further define another cover relation $\looparrowright'$ from $\leq$ as in \autoref{Joins->Covers} above, i.e.
\[C\looparrowright'D\qquad\Leftrightarrow\qquad\forall E,F\in\mathcal{F}(A)\ (\forall c\in C\ (E\leq\{c\}\cup F)\ \Rightarrow\ E\leq D\cup F).\]
We show that the new cover relation $\looparrowright'$ is weaker than the original cover relation $\looparrowright$ in general but agrees $\looparrowright$ whenever the left argument is a singleton.

\begin{prp}
The cover relation $\looparrowright'$ satisfies
\[C\looparrowright D\qquad\Rightarrow\qquad C\looparrowright'D\qquad\Rightarrow\qquad\{c\}\looparrowright D\text{ if }C=\{c\}.\]
\end{prp}

\begin{proof}
Say $C\looparrowright D$ and take any $E,F\in\mathcal{F}(A)$ such that $\{e\}\looparrowright\{c\}\cup F$, for all $e\in E$ and $c\in C$.  Let $C=\{c_1,\cdots,c_n\}$ and note \eqref{Transitive+Monotone} yields $\{e\}\cup C\looparrowright D\cup F$ and $\{e\}\cup C\setminus\{c_1\}\looparrowright\{c_1\}\cup D\cup F$ and hence $\{e\}\cup C\setminus\{c_1\}\looparrowright D\cup F$.  Then \eqref{Transitive+Monotone} again yields $\{e\}\cup C\setminus\{c_1,c_2\}\looparrowright\{c_2\}\cup D\cup F$ and hence $\{e\}\cup C\setminus\{c_1,c_2\}\looparrowright D\cup F$.  Continuing in this way yields $\{e\}\looparrowright D\cup F$, for all $e\in E$, and hence $E\leq D\cup F$, showing that $C\looparrowright D$ implies $C\looparrowright'D$.

Now if $\{c\}\looparrowright'D$ then taking $E=\{c\}$ and $F=\emptyset$ in the definition of $\looparrowright'$, we certainly have $\{c\}\leq\{c\}\cup F$ and hence $\{c\}\leq D$, i.e. $\{c\}\looparrowright D$.
\end{proof}

So the only remaining issue is that a cover relation may not be completely determined by its corresponding $\vee$-semilattice.  However, this problem disappears for cover relations satisfying an appropriate analog of subfitness, namely
\[\tag{Subfit}C\not\looparrowright D\quad\Rightarrow\quad\forall p\in S\ \exists F\ (\forall c\in C\ (\{p\}\looparrowright\{c\}\cup F)\text{ and }C\not\looparrowright D\cup F).\]
Indeed, taking any $p\in D$ here yields $F$ such that $\{p\}\looparrowright\{c\}\cup F$ but $C\not\looparrowright D\cup F$ and hence $\{p\}\not\looparrowright D\cup F$, showing that $C\not\looparrowright D$ implies $C\not\looparrowright'D$ so $\looparrowright\ =\ \looparrowright'$.  It follows that the theory we have developed for semilattice subbases could be translated to arbitrary subbases with their canonical cover structure.

\newpage

\bibliography{maths}{}
\bibliographystyle{alphaurl}

\end{document}